\documentclass[aap,MSNbibl,seceqn,citesort,dvips]{arximspdf}

\doi{10.1214/12-AAP853} 
\volume{22}
\issue{6}
\pubyear{2012}
\firstpage{2539}
\lastpage{2559}

\makeatletter
\newtheorem{theorem}{Theorem}[section]
\newtheorem{lemma}[theorem]{Lemma}
\newtheorem{proposition}[theorem]{Proposition}
\newproclaim{remark}[theorem]{Remark}
\newtheorem{corollary}[theorem]{Corollary}

\newcommand{\F}{\bar{F}}
\newcommand{\Fi}{\bar{\Phi}}
\makeatother

\begin{document}
\begin{frontmatter}

\title{On the rate of approximation in finite-alphabet longest
increasing subsequence problems}
\runtitle{Approximation rate in LIS problems}

\begin{aug}
\author[A]{\fnms{Christian} \snm{Houdr\'e}\thanksref{t1}\ead[label=e1]{houdre@math.gatech.edu}}
\and
\author[B]{\fnms{Zsolt} \snm{Talata}\corref{}\thanksref{t2}\ead[label=e2]{talata@math.ku.edu}}
\thankstext{t1}{Supported in part by NSA Grant H98230-09-1-0017.}
\thankstext{t2}{Supported in part by NSF Grant DMS-09-06929.}
\runauthor{C. Houdr\'e and Z\textup{s}. Talata}
\affiliation{Georgia Institute of Technology and University of Kansas}
\address[A]{School of Mathematics\\
Georgia Institute of Technology\\
Atlanta, Georgia 30332-0160\\
USA\\
\printead{e1}} 
\address[B]{Department of Mathematics\\
University of Kansas\\
405 Snow Hall\\
1460 Jayhawk Boulevard\\
Lawrence, Kansas 66045-7523\\
USA\\
\printead{e2}}
\end{aug}

\received{\smonth{2} \syear{2012}}

%
\begin{abstract}
The rate of convergence of the distribution of the length of the
longest increasing subsequence, toward the maximal eigenvalue of
certain matrix ensembles, is investigated. For finite-alphabet uniform
and nonuniform i.i.d. sources, a rate of $\log n /\sqrt{n}$ is
obtained. The uniform binary case is further explored, and an improved
$1/\sqrt{n}$ rate obtained.
\end{abstract}

%


\end{frontmatter}

\section{Introduction}

In this paper, we consider the length of the longest increasing
subsequence of a random word of size $n$ for general i.i.d. sequences
with alphabet of fixed size $m$. As $n\to\infty$, the limiting
distribution of the normalized length has direct connections to random
matrix theory. If the i.i.d. sequence is uniformly distributed, Tracy
and Widom~\cite{TW} proved that the limiting distribution is that of
the largest eigenvalue of the $m\times m$ traceless Gaussian Unitary
Ensemble (GUE); while for general i.i.d. sequences, Its, Tracy and
Widom~\cite{ITW1,ITW2} showed it to be the distribution of the largest
eigenvalue of a direct sum of certain elements of GUEs.

Limiting distributions in similar problems have also been formulated as
Brownian functionals~\cite{Bar,GW,GTW,LIS}. In particular, in \cite
{LIS}, the length of the longest increasing subsequence is obtained as
a random walk functional, and the limiting distribution, as a Brownian
functional. This direct approach allows us to explore several questions
of probabilistic and statistical nature in this problem, such as the
investigation of the rate of convergence to the limiting distribution,
which is done below.

To briefly describe the content of the paper, for general i.i.d.
sequences we derive, in Section~\ref{secgen}, an upper bound of order
$\log n /\sqrt{n}$ on the rate of convergence, using strong
approximation techniques. In the special case of uniform binary
sequences ($m=2$), the rate is sharpened in Section~\ref{secbin} to
the order $1/\sqrt{n}$.

In previous and related studies, the rate of convergence of certain
random walk functionals has been investigated. For example, in queuing
theory, Glynn and Whitt~\cite{GW} obtained a similar rate via the KMT
technique. In that problem, although the functional of an
$m$-dimensional random walk is similar, the random walks are mutually
independent, which is not our case. Moreover, what is meant there by
rate of convergence is an almost sure upper bound on the deviation
between the random walk and the Brownian functional. The order of that
bound is given, but its constant factor may depend on the realization
of the process. Here, the random walks are dependent, and by rate of
convergence, we mean an upper bound on the deviation of the
distribution functions.

Although the Skorokhod embedding of random walks usually provides a
rate of $\mathcal{O}(n^{-1/4})$~\cite{lif}, when the random walk is one
dimensional and the functionals are the supremum or the local score,
Etienne and Vallois~\cite{EV} obtained a rate of $\mathcal{O}(\sqrt
{\log n/ n})$ using embedding techniques. It is not clear whether or
not their results can be used or generalized to our problem.

To begin we introduce in the next section some notation and also
summarize some of the interplay between the longest increasing
subsequence problem and random matrix theory.

\section{Longest increasing subsequences}\label{secLIn}

Let $X_1, X_2, \ldots$ be a sequence of i.i.d. random variables with
values in the ordered alphabet $\mathcal{A}=\{\alpha_1, \ldots, \alpha
_m\}$, where $\alpha_1<\alpha_2<\cdots<\alpha_m$. Let $p_r=\mathbb
{P}(X_1=\alpha_r)$, $r=1,\ldots,m$, with $p_{\mathrm{max}} = \max_{1\le
r\le m} p_r$ and let also $k$ be the multiplicity of $p_{\mathrm{max}}$
among the probabilities $p_r$ $(1\le r\le m)$.
%
\begin{equation}\label{eqkdef}
k= \# \{ r\dvtx  1\le r\le m, p_r = p_{\mathrm{max}} \}.
\end{equation}
Finally, let $LI_n$ be the length of the longest increasing subsequence
of $X_1, \ldots,\break X_n$, that is,
\[
LI_n = \max\{j\dvtx  X_{i_1} \le X_{i_2} \le\cdots\le X_{i_j},
\mbox{ for some } 1\le i_1 < i_2 <\cdots< i_j \le n \}
.
\]
Properly renormalized, $LI_n$ is known to converge to the maximal
eigenvalue of some matrix ensemble (see~\cite{ITW1,ITW2,Joh,TW}). In fact, in the notation of~\cite{LIS},
%
\begin{equation}\label{eqconv}
\frac{LI_n - n p_{\mathrm{max}} }{\sqrt{n p_{\mathrm{max}}}} \Rightarrow J_k
,
\end{equation}
where
%
\begin{equation}\label{eqH}
\sqrt{p_{\mathrm{max}}}  J_k = -\frac1m \sum_{r=1}^{m-1} r \sigma_r
\tilde{B}^r (1)
+ \mathop{\mathop{\max_{0=t_0\le t_1\le\cdots}}_{\le t_{m-1}\le
t_m=1}}_{
t_r=t_{r-1}, r\in I^*}
\sum_{r=1}^{m-1} \sigma_r \tilde{B}^r (t_r)
,
\end{equation}
with $\sigma_r^2 = p_r + p_{r+1} - (p_r - p_{r+1})^2$, $r=1,2,\ldots
,m-1$ and $I^*=\{ r\dvtx  p_r<p_{\mathrm{max}}, 1\le r\le m \}$. Above,
$(\tilde{B}^1(t),\ldots,\tilde{B}^{m-1}(t))^\top$ is an
$(m-1)$-dimensional driftless Brownian motion with covariance matrix
\[
t \pmatrix{
1 & \rho_{1,2} & \rho_{1,3} & \cdots& \rho_{1,m-1} \vspace*{2pt}\cr
\rho_{2,1} & 1 & \rho_{2,3} & \cdots& \rho_{2,m-1} \vspace*{2pt}\cr
\vdots& \vdots& \ddots& \ddots& \vdots\vspace*{2pt}\cr
\vdots& \vdots&  & 1 & \rho_{m-2,m-1} \vspace*{2pt}\cr
\rho_{m-1,1} & \rho_{m-1,2} & \cdots& \rho_{m-1,m-2} & 1}
,
\]
where
\[
\rho_{r,s} = \cases{
\displaystyle- \frac{p_r + \mu_r\mu_s}{\sigma_r\sigma_s}, & \quad $\mbox{if }
s=r-1,$\vspace*{2pt}\cr
\displaystyle- \frac{p_s + \mu_r\mu_s}{\sigma_r\sigma_s}, & \quad$\mbox{if }
s=r+1,$\vspace*{2pt}\cr
\displaystyle- \frac{\mu_r\mu_s}{\sigma_r\sigma_s}, &\quad $\mbox{if } |r-s|>1, 1\le r,s \le
m-1,$}
\]
and $\mu_r=p_r-p_{r+1}$, $1\le r\le m-1$.

Next, let
%
\begin{equation}\label{eqHtilde}
H_m = \sqrt2 \Biggl\{ -\frac1m \sum_{r=1}^{m-1} r \bar{B}^r (1)
+ \mathop{\max_{0\le t_1\le\cdots}}_{\le t_{m-1}\le1}
\sum_{r=1}^{m-1} \bar{B}^r (t_r) \Biggr\}
,
\end{equation}
where $(\bar{B}^1(t),\ldots,\bar{B}^{m-1}(t))^\top$ is an
$(m-1)$-dimensional driftless Brownian motion with covariance matrix
\[
t \pmatrix{
1 & -1/2 &  &  & \bigcirc\vspace*{2pt}\cr
-1/2 & 1 & -1/2 &  &  \vspace*{2pt}\cr
 & \ddots& \ddots& \ddots&  \vspace*{2pt}\cr
\bigcirc&  & -1/2 & 1 & -1/2 \vspace*{2pt}\cr
 &  &  & -1/2 & 1}
.
\]
Comparing (\ref{eqHtilde}) and (\ref{eqH}), it is immediate that if
the distribution on the alphabet~$\mathcal{A}$ is uniform, that is,
$p_r=1/m$, $r=1,\ldots,m$, then $k=m$, $\mu_r=0$, $\sigma_r^2=2/m$, and
thus $J_m= H_m$, and therefore
\[
\frac{LI_n - n/m}{\sqrt{n/m}} \Rightarrow H_m
.
\]
Actually, similar results hold true for countable
alphabets (see~\cite{LIS}) and our methodology also gives the rate in that case.

Let us now briefly recall the connections, originating in~\cite{Bar}
and~\cite{GTW}, between random matrix theory and the Brownian
functionals encountered in the present paper.

An $m\times m$ element of the Gaussian Unitary Ensemble (GUE) is an
$m\times m$ Hermitian random matrix $\{ Y_{i,j} \}_{ 1\le i,j \le m }$
with $Y_{i,i}\sim N(0,1)$ for $1\le i\le m$, $\mathrm{Re} (Y_{i,j})\sim
N(0,1/2)$ and $\mathrm{Im} (Y_{i,j})\sim N(0,1/2)$ for $1\le i<j\le m$,
and $Y_{i,i}$, $\mathrm{Re} (Y_{i,j})$, $\mathrm{Im} (Y_{i,j})$ are
mutually independent for $1\le i\le j\le m$.

Writing $x^{(m)} = (x_1,x_2,\dots,x_m) \in\mathbb{R}^m$ for any $m \ge
1$, letting $\Delta(x^{(m)}) = \Pi_{1\le i < j \le m}(x_i - x_j)$ be
the Vandermonde determinant, the following facts hold true.

First, from~\cite{TW} and~\cite{LIS}, $\lambda_1^{(m,0)} \stackrel{\mathcal{L}}{=} H_m$, where $\lambda_1^{(m,0)}$ is the largest eigenvalue of
the $m \times m$ \textit{traceless} GUE. Using the joint density of the
eigenvalues of the traceless $m \times m$ GUE~\cite{Mehta,TW}, the distribution function of $H_m$ can be computed directly, for
all $m \ge2$ and all $s \ge0$, as
%
\begin{equation}\label{eqHdens}
\mathbb{P}(H_m \le s) = c_m^0 \int _{\{\max x_j \le s\}}
e^{-({1}/{2})\sum_{i=1}^{m}x_i^2} \Delta\bigl(x^{(m)}\bigr)^2  \lambda_m\bigl(dx^{(m)}\bigr)
,
\end{equation}
where $\lambda_m$ is the Lebesgue measure concentrated on the
hyperplane $\mathcal{L}_m=\{x \in\mathbb{R}^m\dvtx \sum_{i=1}^m x_i = 0\}$,
and where
\[
(c_m^0)^{-1} = \int_{\mathbb{R}^m} e^{-({1}/{2})\sum_{i=1}^{m}x_i^2}
\Delta\bigl(x^{(m)}\bigr)^2  \lambda_m\bigl(dx^{(m)}\bigr)
= (2\pi)^{(m-1)/2}  \prod_{i=0}^{m-1}  i!
.
\]
Note that $H_m$ is a.s. nonnegative, and so $\mathbb{P}(H_m \le s) =
0$ for all $s < 0$.

Second, for all $k\ge2$, $J_k$ can be represented~\cite{LIS,InHom} as
%
\begin{equation}\label{eqsum}
J_k = H_k + \sqrt{ \frac{1-k p_{\mathrm{max}} }{k} }  Z
,
\end{equation}
where $Z$ is a standard normal random variable and, moreover, $H_k$ and
$Z$ are independent, while, $J_1 = \sqrt{ 1- p_{\mathrm{max}} }  Z$.

The distribution of $J_k$ can be described~\cite{ITW1,ITW2}
as the largest eigenvalue of the direct sum of $d$ mutually independent
GUEs, each of size $k_j \times k_j$, $1 \le j \le d$, subject to the
eigenvalue constraint $\sum_{i=1}^m \sqrt{p_i}\lambda_i = 0$. The $k_j$
are the multiplicities of the probabilities having common values, the
$p_i$ are ordered in nonincreasing order and the eigenvalues are
ordered in terms of the GUEs corresponding to the appropriate values of $p_i$.

As shown in~\cite{ITW2}, for any $k \ge1$ and all $s \in\mathbb{R}$,
$J_k$ has distribution given by
\begin{eqnarray*}
&&\mathbb{P}(J_k \le s)\\
&&\quad = c_{k,p_{\mathrm{max}}} \int _{\{\max x_j \le s\}}
e^{-({1}/{2})[\sum_{i=1}^{k}x_i^2
+ ({p_{\mathrm{max}}}/{(1-kp_{\mathrm{max}})})(\sum_{i=1}^{k}x_i)^2]} \Delta
\bigl(x^{(k)}\bigr)^2  \,dx^{(k)}
,
\end{eqnarray*}
where
\[
c_{k,p_{\mathrm{max}}}^{-1} = \int_{\mathbb{R}^k}
e^{-({1}/{2})[\sum_{i=1}^{k}x_i^2
+ ({p_{\mathrm{max}}}/({1-kp_{\mathrm{max}}}))(\sum_{i=1}^{k}x_i)^2]}
\Delta\bigl(x^{(k)}\bigr)^2  \,dx^{(k)}
.
\]

Below, we study the rate of approximation in (\ref{eqconv}) and prove
(see Section~\ref{secgen}) that
\[
\sup_{x\in\mathbb{R}}
\biggl|
\mathbb{P}( \frac{LI_n - n p_{\mathrm{max}}}{\sqrt{n p_{\mathrm{max}}}} \ge x )
- \mathbb{P}( J_k \ge x )
 \biggr|
\le C(m,k)  \frac{\log n}{\sqrt{n}}
,
\]
where the constant $C(m,k)$ depends only on $m$ and $k$.

\section{Upper bounds on the density functions}\label{secdens}

Our first results provide upper bounds on the density of the
functionals $J_k$ and $H_k$.

\begin{proposition}\label{corJpdf}
\textup{(i)} Let $f_{H_k}$ be the probability density function of $H_k$. Then,
for any $k=2,3,\ldots,m$,
\[
\sup_{x\in\mathbb{R}}  f_{H_k} (x) \le
k^{3k} ( 2\pi e^2 )^{k/2} \sqrt{\frac{e}{\pi}}
.
\]
\begin{longlist}[(ii)]
\item[(ii)] Let $f_{J_k}$ be the probability density function of $J_k$. Then,
for any $k=2,3,\ldots,m-1$,\vspace*{-2pt}
\[
\sup_{x\in\mathbb{R}}  f_{J_k} (x) \le
\min\Biggl\{
\sqrt{\frac{k}{ 2\pi(1-kp_{\mathrm{max}}) }} , k^{3k} ( 2\pi e^2 )^{k/2} \sqrt{\frac{e}{\pi}}
\Biggr\}\vspace*{-1pt}
\]
and for $k=1$, $\sup_{x\in\mathbb{R}}  f_{J_1} (x) = 1/\sqrt{ 2\pi
(1-p_{\mathrm{max}}) }$.
\end{longlist}
\end{proposition}

\begin{remark}
The distribution function of the largest eigenvalue of the $k\times k$
GUE can be computed directly~\cite{Mehta}, and so does the one of the
$k\times k$ traceless GUE in (\ref{eqHdens}). Its derivative, the
density function, is upper bounded by $k$ times the density function of
the one-dimensional marginal of the distribution of the eigenvalues of
the $k\times k$ GUE. Both the joint density of the eigenvalues of the
$k\times k$ GUE and its marginals have a determinantal representation
using Hermite polynomials~\cite{Mehta}, Section 6.2, which
provides an upper bound of order $k^{5/6}$ on the density of the largest
eigenvalue of the $k\times k$ GUE. In turn, this gives the order of the constant with the rate of
convergence $\log n/\sqrt{n}$, in the framework of~\cite{GW}.

For the joint and marginal density of the eigenvalues of the $k\times
k$ traceless GUE, a determinantal representation does not seem to be
available. One can still conjecture a polynomial upper bound on the
supremum of the density of the largest eigenvalue of the $k\times k$
traceless GUE, but the authors' efforts did not lead to such a bound in
part (i) of Proposition~\ref{corJpdf}. Indeed, the traceless
condition induces dependencies between the entries of the Gaussian
random matrix making the analysis more delicate than in the GUE
case.\looseness=1
\end{remark}

\begin{pf} First for (ii) using (\ref{eqsum}), for $k>1$,\vspace*{-1pt}
\begin{eqnarray*}
f_{J_k}(x)& =& \int_{\mathbb{R}} f_{H_k} (u)  f_{\sqrt{(1-k p_{\mathrm{max}}) /k } Z} (x-u)  \,du
\\[-2pt]
&\le&\sup_{u\in\mathbb{R}} f_{H_k} (u)
,\vspace*{-1pt}
\end{eqnarray*}
and use (i). Similarly, for $k<m$,\vspace*{-1pt}
\[
\sup_{x\in\mathbb{R}}  f_{J_k} (x)
\le\sup_{u\in\mathbb{R}} f_{\sqrt{(1-k p_{\mathrm{max}}) /k } Z} (u)
= \sqrt{\frac{k}{ 2\pi(1-kp_{\mathrm{max}}) }}
.\vspace*{-1pt}
\]
Now for (i), to upper bound the density of $H_k$, consider its
cumulative distribution function. By (\ref{eqHdens}),\vspace*{-1pt}
\begin{eqnarray*}
\mathbb{P}(H_k \le s)
&=& \frac{ 1 }{ (2\pi)^{(k-1)/2} }
\\[-3pt]
&&{}\times\int _{\{\max x_j \le s\}} \exp\Biggl( -\frac12 \sum
_{i=1}^k x_i^2 +2 \sum _{ 1\le i<j \le k } \log|x_i-x_j| \\[-3pt]
&&\hspace*{165pt}\qquad{}- \sum
_{i=1}^{k-1} \log i! \Biggr) \lambda_k(dx)
,\vspace*{-1pt}
\end{eqnarray*}
and so
\begin{eqnarray}\label{eqcdf2}
&&\mathbb{P}(s\le H_k \le s+\varepsilon)
\nonumber\\
&&\qquad= \int _{\{s\le\max x_j \le s+\varepsilon\}}
\frac{ e^{-(1/2) ( 1-{2}/{\nu^2} ) \sum_{i=1}^k
x_i^2} }{ (2\pi)^{(k-1)/2} }
\\
&&\qquad\quad{}\times \exp\Biggl( -\frac{1}{\nu^2} \sum_{i=1}^k x_i^2 +2 \sum
_{ 1\le i<j \le k } \log|x_i-x_j| - \sum_{i=1}^{k-1} \log i! \Biggr)
 \lambda_k(dx)
\nonumber
\end{eqnarray}
for any $\nu>\sqrt{2}$. In order to dominate the first term of the
integrand in (\ref{eqcdf2}), a~bound (see~\cite{Mehta}, Appendix~A.6) going
back to Stieltjes, asserts that
\begin{eqnarray*}
&&\frac12 \sum_{i=1}^k x_i^2 - \sum_{ 1\le i<j \le k } \log|x_i-x_j|\\
&&\qquad \ge
\frac14 k(k-1)(1+\log2) - \frac12 \sum_{i=1}^k i \log i
.
\end{eqnarray*}
Hence,
\begin{eqnarray}\label{eqcdfb1}
&&\frac{1}{\nu^2} \sum_{i=1}^k x_i^2 - 2 \sum_{ 1\le i<j \le k } \log|x_i-x_j|
\nonumber\\
&&\qquad= 2\Biggl( \frac{1}{2\nu^2} \sum_{i=1}^k x_i^2 - \sum_{ 1\le i<j \le k
} \log|x_i-x_j| \Biggr)
\nonumber\\
&&\qquad= 2\Biggl( \frac12 \sum_{i=1}^k \biggl(\frac{x_i}{\nu}\biggr)^2 - \sum
_{ 1\le i<j \le k } \log\biggl| \frac{x_i}{\nu} - \frac{x_j}{\nu}
\biggr| \Biggr) - \frac{k(k-1)}{2}  2\log\nu
\\
&&\qquad\ge2\Biggl( \frac14 k(k-1)(1+\log2) - \frac12 \sum_{i=1}^k i \log i
\Biggr) - k(k-1)\log\nu
\nonumber
\\
&&\qquad= \frac12 k(k-1)(1+\log2) - \sum_{i=1}^k i \log i - k(k-1)\log\nu.
\nonumber
\end{eqnarray}
On the other hand,
%
\begin{eqnarray}\label{eqcdfb2}
\qquad\sum_{i=1}^{k-1} \log i! &= &\sum_{i=1}^{k-1} \sum_{j=1}^i \log j = \sum
_{j=1}^{k-1} (k-j) \log j
\nonumber
\\[-8pt]
\\[-8pt]
\nonumber
&=& k \log(k-1)! - \sum_{i=1}^{k-1} i \log i.
\end{eqnarray}
Combining (\ref{eqcdfb1}) and (\ref{eqcdfb2}) leads to
\begin{eqnarray}\label{eqcdffirst}
&&\exp\Biggl( -\frac{1}{\nu^2} \sum_{i=1}^k x_i^2 +2 \sum _{ 1\le
i<j \le k } \log|x_i-x_j| - \sum_{i=1}^{k-1} \log i! \Biggr)
\nonumber\\
&&\qquad\le\exp\Biggl( - \frac12 k(k-1)(1+\log2) + k(k-1) \log\nu- k \log(k-1)!
\vphantom{\sum_{i=1}^{k-1} i \log i +k\log k}\\
& &\hspace*{175pt}\qquad{}+ 2 \sum_{i=1}^{k-1} i \log i +k\log k \Biggr).
\nonumber
\end{eqnarray}
Here, using Stirling's inequality (see Feller~\cite{Feller}, page~54)
\[
(k-1)! > \sqrt{2\pi(k-1)}  (k-1)^{(k-1)} e^{-(k-1)+1/(12(k-1)+1)}
,
\]
and
\[
\sum_{i=1}^{k-1} i \log i \le\int_1^k i \log i \, di
= \frac{k^2}{2} \log k - \frac{k^2}{4} + \frac14
,
\]
the exponent in (\ref{eqcdffirst}) can be upper bounded by
\begin{eqnarray*}
&&- \frac12 k(k-1)(1+\log2) + k(k-1)\log\nu- k(k-1) \log(k-1)
\\
&&\quad{}-k \biggl( -(k-1)+\frac{1}{12(k-1)+1} + \log\sqrt{2\pi(k-1)} \biggr)
+k\log k
\\
&&\quad{}+ k^2 \log k - \frac{k^2}{2} + \frac12
\\
&&\qquad= k^2 \biggl( -\frac12 (1+\log2) + \log\nu- \log(k-1) + 1 + \log k -
\frac12 \biggr)
\\
&&\qquad\quad{}+ k \biggl( \frac12 (1+\log2) - \log\nu+ \log(k-1)\\
&&\hspace*{23pt}\qquad\quad{}
- 1 - \frac{1}{12(k-1)+1}
+ \log\sqrt{2\pi(k-1)} + \log k \biggr) +\frac12
\\
&&\qquad\le k^2 \biggl( \log\frac{\nu}{\sqrt{2}} + \log\frac{k}{k-1} \biggr)
+ \frac52 k\log k + k \log\sqrt{\frac{4\pi}{\nu^2 e}} + \frac12
.
\end{eqnarray*}
And the bound (\ref{eqcdffirst}) becomes
\begin{eqnarray}\label{eqcdfsec}
&&\exp\Biggl( -\frac14 \sum_{i=1}^k x_i^2 +2 \sum _{ 1\le i<j \le
k } \log|x_i-x_j| - \sum_{i=1}^{k-1} \log i! \Biggr)\nonumber\\
&&\qquad\le\biggl(\frac{\nu}{\sqrt{2}}\biggr)^{k^2} \biggl( \frac{k}{k-1}
\biggr)^{k^2} k^{5k/2} \biggl( \frac{4\pi}{\nu^2 e} \biggr)^{k/2} \sqrt{e}\\
&&\qquad\le\biggl(\frac{\nu}{\sqrt{2}}\biggr)^{k^2} k^{5k/2} \biggl( \frac
{4\pi e}{\nu^2} \biggr)^{k/2} \sqrt{e}.
\nonumber
\end{eqnarray}
Thus (\ref{eqcdf2}) becomes
\begin{eqnarray}\label{eqcdf3}
\qquad\mathbb{P}(s\le H_k \le s+\varepsilon)&\le&\biggl(\frac{\nu}{\sqrt{2}}\biggr)^{k^2} k^{5k/2} \biggl( \frac
{4\pi e}{\nu^2} \biggr)^{k/2}
\nonumber
\\[-8pt]
\\[-8pt]
\nonumber
&&{}\times\sqrt{e}
\int _{\{s\le\max x_j \le s+\varepsilon\}} \frac{ e^{-(1/2)
( 1-{2}/{\nu^2} ) \sum_{i=1}^k x_i^2} }{ (2\pi
)^{(k-1)/2} }  \lambda_k(dx)
.
\end{eqnarray}
On the hyperplane $\mathcal{L}_k \subset\mathbb{R}^k$, the function
\[
\frac{ 1 }{ (2\pi)^{(k-1)/2}  [(1-2/\nu^2)^{-1}]^{(k-1)/2} }
e^{-(1/2) ( 1-{2}/{\nu^2} ) \sum_{i=1}^k x_i^2}
\]
is the probability density function of the $(k-1)$-dimensional normal
distribution with mean $(0,\ldots,0)\in\mathbb{R}^{k-1}$ and covariance
$[(1-2/\nu^2)^{-1}] I_{k-1}$, where $I_{k-1}$ is the
$(k-1)$-dimensional identity matrix. Therefore,
\begin{eqnarray}\label{eqcdf4}
&&\int _{\{s\le\max x_j \le s+\varepsilon\}} \frac{ e^{-(1/2)( 1-{2}/{\nu^2} ) \sum_{i=1}^k x_i^2} }{ (2\pi
)^{(k-1)/2} }  \lambda_k(dx)
\nonumber\\
&&\qquad= \biggl( 1-\frac{2}{\nu^2} \biggr)^{-(k-1)/2} \nonumber\\
&&\qquad\quad{}\times\int _{\{s\le
\max x_j \le s+\varepsilon\}} \frac{ e^{-(1/2) ( 1-{2}/{\nu
^2} ) \sum_{i=1}^k x_i^2} }{ (2\pi)^{(k-1)/2}  [(1-2/\nu
^2)^{-1}]^{(k-1)/2} }  \lambda_k(dx)
\nonumber\\
&&\qquad\le\biggl( 1-\frac{2}{\nu^2} \biggr)^{-(k-1)/2} \\
&&\qquad\quad{}\times\sum_{j=1}^{k} \int
 _{\{s\le x_j \le s+\varepsilon\}} \frac{ e^{-(1/2)
(1-{2}/{\nu^2} ) \sum_{i=1}^k x_i^2} }{ (2\pi)^{(k-1)/2}
[(1-2/\nu^2)^{-1}]^{(k-1)/2} }  \lambda_k(dx)
\nonumber\\
&&\qquad\le\biggl( 1-\frac{2}{\nu^2} \biggr)^{-(k-1)/2} \sum_{j=1}^{k} \sqrt
{2}  \varepsilon\sup _{\{s\le x_j \le s+\varepsilon\}} \frac{
e^{-(1/2) ( 1-{2}/{\nu^2}) x_j^2 } }{ \sqrt{ 2\pi
(1-2/\nu^2)^{-1} } }
\nonumber
\\
&&\qquad\le \frac{\varepsilon}{\sqrt{\pi}}  k \biggl( 1-\frac{2}{\nu^2}
\biggr)^{-k/2+1}
\nonumber
.
\end{eqnarray}
Using (\ref{eqcdf4}), (\ref{eqcdf3}) yields
\begin{eqnarray*}\label{eqcdf5}
&&\mathbb{P}(s\le H_k \le s+\varepsilon)
\\
&&\qquad\le\biggl(\frac{\nu^2}{2}\biggr)^{k^2/2} \biggl( 1-\frac{2}{\nu^2}
\biggr)^{-k/2+1}
\varepsilon k^{5k/2+1} \biggl( \frac{4\pi e}{\nu^2} \biggr)^{k/2} \sqrt
{\frac{e}{\pi}}.
\end{eqnarray*}
Choosing $\nu>\sqrt{2}$ such that $1-2/\nu^2=1/k$ provides
\[
\biggl(\frac{\nu^2}{2}\biggr)^{k^2/2} = \biggl(1-\frac{1}{k}
\biggr)^{-k^2/2} \le e^{k/2}
,
\]
leading to
\[
\mathbb{P}(s\le H_k \le s+\varepsilon)
\le\varepsilon k^{3k} ( 2\pi e^2 )^{k/2} \sqrt{\frac{e}{\pi}}
,
\]
and the proof is complete.
\end{pf}

\section{Rate of convergence results}\label{secgen}

Below, we study the rate of convergence in~(\ref{eqconv}) and show that:

\begin{theorem}\label{thgen}
For any $n\in\mathbb{N}$, $m\in\mathbb{N}$, for $k=2,3,\ldots,m$,
\begin{eqnarray}\label{eqTh}
&&\sup_{x\in\mathbb{R}}
\biggl|
\mathbb{P}\biggl( \frac{LI_n - n p_{\mathrm{max}}}{\sqrt{n p_{\mathrm{max}}}} \ge x \biggr)
- \mathbb{P}( J_k \ge x )
 \biggr|
\nonumber\hspace*{-35pt}
\\[-4pt]
\\[-12pt]
\nonumber
&&\qquad\le c (m-1) \Biggl( (m-1)^2 \sigma_{\mathrm{max}}^2 + \Biggl( \frac
{k^{3k}}{\sqrt{p_{\mathrm{max}}}} \wedge\sqrt{\frac{k}{ p_{\mathrm{max}}
(1-kp_{\mathrm{max}}) }} \Biggr) \Biggr) \frac{\log n}{\sqrt{n}},\hspace*{-35pt}
\end{eqnarray}
where $c>0$ is an absolute constant and $\sigma_{\mathrm{max}} = \max
_{1\le r\le m-1} \sigma_r$.
For $k=1$, (\ref{eqTh}) holds with the minimum replaced by $(p_{\mathrm{max}} (1-p_{\mathrm{max}}))^{-1/2}$.
\end{theorem}

\begin{remark}
For $k=m$, the minimum in (4.14) is $m^{3m+1/2}$ since\break
$\sqrt{k / ( p_{\mathrm{max}} (1-kp_{\mathrm{max}}) )}$ is then
understood to be infinite. In general, the minimum is $k^{3k} /\sqrt
{p_{\mathrm{max}}}$, if $p_{\mathrm{max}} \ge(k^{6k} - k)/(k^{6k+1})$ and
$\sqrt{k / ( p_{\mathrm{max}} (1-kp_{\mathrm{max}}) )}$, if
$p_{\mathrm{max}} \le(k^{6k} - k)/(k^{6k+1})$.
\end{remark}

In particular, (\ref{eqTh}) implies the following result which should
be contrasted with Theorem 4 and Theorem 6 of~\cite{BH}.

\begin{corollary}\label{corgenconv}
If $k$ is fixed and $m\to\infty$ as $n\to\infty$ in such a way that
$m=o(n^{1/4} \log^{-1/2} n)$, then
\[
\frac{LI_n - n p_{\mathrm{max}}}{\sqrt{n p_{\mathrm{max}}}} \Rightarrow J_k
.
\]
\end{corollary}

\begin{pf}
This immediately follows from Theorem~\ref{thgen} since $\sigma_{\mathrm{max}}^2 \le2/m$ and $1/m \le p_{\mathrm{max}} \le1/k$.
\end{pf}

\begin{pf*}{Proof of Theorem \protect\ref{thgen}}
Set $L_n = (LI_n - n p_{\mathrm{max}})/\sqrt{n}$, and for $i=1,\ldots,n$
and $r=1,\ldots,m-1$, set also
\[
Z^r_i= \cases{
1, & \quad$\mbox{if } X_i=\alpha_r,$
\vspace*{2pt}\cr
-1, & \quad$\mbox{if } X_i=\alpha_{r+1},$
\vspace*{2pt}\cr
0, & \quad$\mbox{otherwise}.$}
\]
Clearly, $\operatorname{Var}  Z_i^r = \sigma_r^2$ and $\mathbb{E} Z_i^r =\mu
_r$. Set $\tilde{S}^r_0=0$ and
\[
\tilde{S}^r_j = \sum_{i=1}^j \frac{ Z^r_i - \mu_r }{ \sigma_r }
,\qquad  j=1,\ldots,n
,
\]
and then $L_n$ can be written (see the proof of Theorem 3.1 in \cite
{LIS}), as
%
\begin{equation}\label{eqLn}
L_n = -\frac1m \sum_{r=1}^{m-1} r \sigma_r \frac{ \tilde{S}^r_n }{
\sqrt{n} }
+ \mathop{\mathop{\max_{0=j_0\le j_1\le\cdots}}_{\le j_{m-1}\le
j_m=n}}_{j_r=j_{r-1}, r\in I^*}
\sum_{r=1}^{m-1} \sigma_r \frac{ \tilde{S}^r_{j_r} }{ \sqrt{n} }
+E_n
,
\end{equation}
where for the remainder term $E_n$ we have for any $\varepsilon>0$,
%
\begin{equation}\label{eqEn}
\mathbb{P}(|E_n|\ge\varepsilon) < \varepsilon\bigl( 1+ (m-1)^2 \sigma
_{\mathrm{max}}^2 \bigr)
.
\end{equation}

Letting
%
\begin{equation}\label{eqHnbar}
\widetilde{H}_{n,k} =
-\frac1m \sum_{r=1}^{m-1} r \sigma_r \tilde{B}^r (1)
+ \mathop{\mathop{\max_{0=j_0\le j_1\le\cdots}}_{\le j_{m-1}\le
j_m=n}}_{j_r=j_{r-1}, r\in I^*}
\sum_{r=1}^{m-1} \sigma_r \tilde{B}^r \biggl( \frac{j_r}{n} \biggr)
,
\end{equation}
we have for any $\varepsilon>0$
\begin{eqnarray}\label{eqdec}\qquad
&&\mathbb{P}\bigl( | L_n/\sqrt{p_{\mathrm{max}}} - J_k | >
2\varepsilon\bigr)
\nonumber
\\[-8pt]
\\[-8pt]
\nonumber
&&\qquad\le\mathbb{P}\bigl( | L_n - \widetilde{H}_{n,k} | >
\varepsilon\sqrt{p_{\mathrm{max}}} \bigr) + \mathbb{P}\bigl( |
\widetilde{H}_{n,k} - \sqrt{p_{\mathrm{max}}} J_k | > \varepsilon
\sqrt{p_{\mathrm{max}}} \bigr).
\end{eqnarray}
Using Lemmas~\ref{lemkmt} and~\ref{lemsup} below, (\ref{eqdec}) can
be upper bounded by
\begin{eqnarray}\label{eqplus}
&&\frac{1+ (m-1)^2 \sigma_{\mathrm{max}}^2}{2} \varepsilon\sqrt{p_{\mathrm{max}}} +
\exp\biggl( -\frac{\xi\varepsilon\sqrt{n p_{\mathrm{max}}} }{16 (m-1)}
\biggr) \sum_{r=1}^{m-1} \biggl( 1+ \frac{\sigma_r\sqrt{n}}{1-|\mu_r|}
\biggr)
\nonumber
\\[-8pt]
\\[-8pt]
\nonumber
&&\qquad{}+4 (m-1) n \exp\biggl( \frac{ -\varepsilon^2 n p_{\mathrm{max}} }{8
\sigma_{\mathrm{max}}^2 (m-1)^2} \biggr)
.
\end{eqnarray}
Now, from Proposition~\ref{corJpdf},
\begin{eqnarray}\label{eqfinal}\quad
&&\bigl|  \mathbb{P}\bigl( L_n/\sqrt{p_{\mathrm{max}}} \ge x \bigr) -
\mathbb{P}( J_k \ge x )  \bigr|
\nonumber\\
&&\qquad\le\mathbb{P}\bigl( \bigl| L_n/\sqrt{p_{\mathrm{max}}} -J_k \bigr| \ge
2\varepsilon\bigr)
+ \mathbb{P}( x-2\varepsilon\le J_k \le x+2\varepsilon)
\nonumber
\\
&&\qquad\le\frac{1+ (m-1)^2 \sigma_{\mathrm{max}}^2}{2} \varepsilon\sqrt
{p_{\mathrm{max}}}
+
\exp\biggl( -\frac{\xi\varepsilon\sqrt{n p_{\mathrm{max}}} }{16 (m-1)}
\biggr) \sum_{r=1}^{m-1} \biggl( 1+ \frac{\sigma_r\sqrt{n}}{1-|\mu_r|}
\biggr)
\\
&&\qquad\quad{} +4 (m-1) n \exp\biggl( \frac{ -\varepsilon^2 n p_{\mathrm{max}}
}{8 \sigma_{\mathrm{max}}^2 (m-1)^2} \biggr)
\nonumber
\\
&&\qquad\quad{} +4\varepsilon  \min\Biggl\{
\sqrt{\frac{k}{ 2\pi(1-kp_{\mathrm{max}}) }} , k^{3k} ( 2\pi e^2 )^{k/2} \sqrt{\frac{e}{\pi}}
\Biggr\}.\nonumber
\end{eqnarray}
With
\[
\varepsilon= \frac{16(m-1)}{\xi\sqrt{p_{\mathrm{max}}}}  \frac{\log
n}{\sqrt{n}}
,
\]
the right-hand side of (\ref{eqfinal}) becomes
\begin{eqnarray*}
&&\frac{\log n}{\sqrt{n}} \Biggl(
\frac{ 8(m-1) (1+ (m-1)^2 \sigma_{\mathrm{max}}^2) }{ \xi}
+\frac{1}{\log n} \sum_{r=1}^{m-1}\biggl ( \frac{1}{\sqrt{n}} + \frac
{\sigma_r}{1-|\mu_r|} \biggr)
\\
&&\hspace*{8pt}\qquad{}
+4(m-1) n^{-( 32/ \sigma_{\mathrm{max}}^2 \xi^2) \log n - {\log\log
n}/{\log n}+ 3/2 }
\\
&&\hspace*{24pt}\qquad{}
+ \frac{64(m-1)}{\xi\sqrt{p_{\mathrm{max}}}}
\min\Biggl\{
\sqrt{\frac{k}{ 2\pi(1-kp_{\mathrm{max}}) }} , k^{3k} ( 2\pi e^2 )^{k/2} \sqrt{\frac{e}{\pi}}
\Biggr\}
\Biggr)
,
\end{eqnarray*}
which yields the claim of the theorem.
\end{pf*}

\begin{lemma}\label{lemkmt}
For any $\varepsilon>0$,
\begin{eqnarray*}
&&\mathbb{P}( | L_n - \widetilde{H}_{n,k} | >
\varepsilon)
\\
&&\qquad\le\bigl( 1+ (m-1)^2 \sigma_{\mathrm{max}}^2 \bigr) \frac{\varepsilon
}{2} +
\exp\biggl( -\frac{\xi\varepsilon\sqrt{n}}{16 (m-1)} \biggr) \sum
_{r=1}^{m-1} \biggl( 1+ \frac{\sigma_r\sqrt{n}}{1-|\mu_r|} \biggr)
,
\end{eqnarray*}
where $\xi>0$ is an absolute constant.
\end{lemma}

\begin{pf}
Comparing (\ref{eqLn}) to (\ref{eqHnbar}),
\begin{eqnarray}\label{eqmax}
&&| L_n - \widetilde{H}_{n,k} |
\nonumber\\
&&\qquad\le |E_n| + \sum_{r=1}^{m-1} \biggl( 1+\frac{r}{m} \biggr) \sigma_r
\max_{0\le j\le n} \biggl| \frac{ \tilde{S}^r_j }{ \sqrt{n} } - \tilde
{B}^r\biggl(\frac{j}{n}\biggr) \biggr|
\\
&&\qquad\le |E_n| + 2 \sum_{r=1}^{m-1} \sigma_r \max_{0\le j\le n} \biggl|
\frac{ \tilde{S}^r_j }{ \sqrt{n} } - \tilde{B}^r\biggl(\frac{j}{n}
\biggr)\biggr |\nonumber
.
\end{eqnarray}
For any $\delta>0$ and $0\le r\le m-1$,
%
\begin{equation}\label{eqmax2}\quad
\mathbb{P}\biggl( \max_{0\le j\le n} \biggl| \frac{ \tilde{S}^r_j }{
\sqrt{n} } - \tilde{B}^r\biggl(\frac{j}{n}\biggr) \biggr| > \delta\biggr)
= \mathbb{P}\Bigl( \max_{0\le j\le n} | \tilde{S}^r_j - \tilde
{B}^r(j) | > \delta\sqrt{n} \Bigr).
\end{equation}
Applying Sakhanenko's version of the KMT inequality~\cite{kmt1,kmt2,lif,Sak} to the partial sums $\tilde{S}^r_j$,
$j=0,\ldots,n$, of the i.i.d. random variables $(Z_i^r-\mu_r) / \sigma
_r$, $i=1,\ldots,n$, (\ref{eqmax2}) can be upper bounded by
%
\begin{equation}\label{eqkmt}
\bigl( 1+C_2\sqrt{n} \bigr) \exp\bigl( -C_1 \delta\sqrt{n} \bigr).
\end{equation}
Above, $C_1=\xi\lambda_r$ and $C_2=\lambda_r$, where $\xi$ is an
absolute constant and
\begin{eqnarray*}
\lambda_r
&=& \sup\biggl\{ \lambda\dvtx  \lambda \mathbb{E} \biggl( \biggl|\frac
{Z_i^r-\mu_r}{\sigma_r}\biggr|^3 \exp\biggl\{ \lambda\biggl|\frac
{Z_i^r-\mu_r}{\sigma_r}\biggr| \biggr\} \biggr) \le\mathbb{E} \biggl(
\frac{Z_i^r-\mu_r}{\sigma_r} \biggr)^2 \biggr\}
\\
&=& \sigma_r \sup\bigl\{ \lambda\dvtx  \lambda \mathbb{E} (
|Z_i^r-\mu_r|^3 \exp\{ \lambda|Z_i^r-\mu_r|
\} ) \le\operatorname{Var}  Z_i^r \bigr\}
.
\end{eqnarray*}
Since $|Z_i^r-\mu_r|\le2$, choosing $\lambda=1/4$ gives
\begin{eqnarray*}
&&\frac{1}{4} \mathbb{E} \biggl( |Z_i^r-\mu_r|^3 \exp\biggl\{
\frac{1}{4} |Z_i^r-\mu_r| \biggr\} \biggr)
\\
&&\qquad=\mathbb{E}\biggl ( ( Z_i^r-\mu_r )^2 \frac{|Z_i^r-\mu
_r|}{4}  \exp\biggl\{ \frac{|Z_i^r-\mu_r|}{4} \biggr\}
\biggr)
\\
&&\qquad\le\mathbb{E} \biggl( ( Z_i^r-\mu_r )^2 \frac{1}{2}  \exp
\biggl\{ \frac{1}{2} \biggr\} \biggr)
\\
&&\qquad\le\operatorname{Var}  Z_i^r
,
\end{eqnarray*}
which implies that $\lambda_r \ge\sigma_r/4$. Next, for any $\lambda>
1/\min\{1-\mu_r,1+\mu_r\}$,
\begin{eqnarray*}
\operatorname{Var}  Z_i^r
&=& p_r (1-\mu_r)^2 + p_{r+1} (1+\mu_r)^2
\\
&<& \lambda\min\{|1-\mu_r|,|1+\mu_r|\} \bigl( p_r (1-\mu_r)^2 +
p_{r+1} (1+\mu_r)^2 \bigr)
\\
&\le&\lambda( p_r |1-\mu_r|^3 + p_{r+1} |1+\mu_r|^3 )
\\
&\le&\lambda \mathbb{E} ( |Z_i^r-\mu_r|^3 )
\\
&\le&\lambda \mathbb{E} ( |Z_i^r-\mu_r|^3 \exp\{
\lambda|Z_i^r-\mu_r| \} )
,
\end{eqnarray*}
which implies that $\lambda_r \le\sigma_r/\min\{1-\mu_r,1+\mu_r\} =
\sigma_r/(1-|\mu_r|)$.
Thus, the upper bound (\ref{eqkmt}) becomes
%
\begin{equation}
\biggl( 1+ \frac{\sigma_r\sqrt{n}}{1-|\mu_r|} \biggr) \exp\biggl( -\frac
{\xi\sigma_r}{4}  \delta\sqrt{n} \biggr)
\label{eqkmt2}
.
\end{equation}
Combining (\ref{eqkmt2}) and (\ref{eqEn}) with (\ref{eqmax}),
\begin{eqnarray*}
&&\hspace*{-3pt}\mathbb{P}( | L_n - \widetilde{H}_{n,k} | >
\varepsilon)
\\
&&\hspace*{-3pt}\qquad\le\mathbb{P}\Biggl( |E_n| + 2 \sum_{r=1}^{m-1} \sigma_r \max_{0\le
j\le n} \biggl| \frac{ \tilde{S}^r_j }{ \sqrt{n} } - \tilde{B}^r
\biggl(\frac{j}{n}\biggr) \biggr| > \varepsilon\Biggr)
\\
&&\hspace*{-3pt}\qquad\le\mathbb{P}\biggl( |E_n| > \frac{\varepsilon}{2} \biggr) +
\sum_{r=1}^{m-1}  \mathbb{P}\Biggl( \max_{0\le j\le n} \biggl| \frac{
\tilde{S}^r_j }{ \sqrt{n} } - \tilde{B}^r\biggl(\frac{j}{n}\biggr)
\biggr| > \frac{\varepsilon}{4 \sigma_r (m-1)} \Biggr)
\\
&&\hspace*{-3pt}\qquad\le\frac{\varepsilon}{2} \bigl( 1+ (m-1)^2 \sigma_{\mathrm{max}}^2
\bigr) +
\sum_{r=1}^{m-1}  \biggl( 1+ \frac{\sigma_r\sqrt{n}}{1-|\mu_r|}
\biggr) \exp\biggl( -\frac{\xi\sigma_r}{4}  \frac{\varepsilon}{4 \sigma_r
(m-1)} \sqrt{n} \biggr)
,
\end{eqnarray*}
and the proof is complete.
\end{pf}

\begin{lemma}\label{lemsup}
For any $\varepsilon>0$,
\[
\mathbb{P}\bigl( \bigl| \widetilde{H}_{n,k} - \sqrt{p_{\mathrm{max}}}
J_k\bigr | > \varepsilon\bigr)
\le4 (m-1) n \exp\biggl( \frac{ -\varepsilon^2 n}{8 \sigma_{\mathrm{max}}^2 (m-1)^2} \biggr)
.
\]
\end{lemma}

\begin{pf}
Comparing (\ref{eqH}) and (\ref{eqHnbar}),
\begin{eqnarray}\label{eqdiff}
&&\bigl| \widetilde{H}_{n,k} - \sqrt{p_{\mathrm{max}}} J_k \bigr|
\nonumber\\
&&\qquad\le\sum_{r=1}^{m-1} \biggl( 1+\frac{r}{m} \biggr) \sigma_r \max_{0\le
j\le n-1}  \sup_{0\le t\le1/n} \biggl| \tilde{B}^r\biggl(\frac
{j}{n}+t\biggr) - \tilde{B}^r\biggl(\frac{j}{n}\biggr) \biggr|
\\
&&\qquad\le2 \sum_{r=1}^{m-1} \sigma_r \max_{0\le j\le n-1}  \sup_{0\le t
\le1/n} \biggl| \tilde{B}^r\biggl(\frac{j}{n}+t\biggr) - \tilde{B}^r
\biggl(\frac{j}{n}\biggr) \biggr|
\nonumber
.
\end{eqnarray}
Here, for any $\delta>0$ and $0\le r\le m-1$,
\begin{eqnarray}\label{equpper}
&&\mathbb{P}\biggl( \max_{0\le j\le n-1}  \max_{0\le t\le1/n} \biggl|
\tilde{B}^r\biggl(\frac{j}{n}+t\biggr) - \tilde{B}^r\biggl(\frac
{j}{n}\biggr) \biggr| > \delta\biggr)
\nonumber\\
&&\qquad\le\sum_{j=0}^{n-1}  \mathbb{P}\biggl( \max_{0\le t\le1/n} \biggl|
\tilde{B}^r\biggl(\frac{j}{n}+t\biggr) - \tilde{B}^r\biggl(\frac
{j}{n}\biggr) \biggr| > \delta\biggr)
\nonumber
\\
&&\qquad= \sum_{j=0}^{n-1}  \mathbb{P}\Bigl( \max_{0\le t\le1} | \tilde
{B}^r(t) | > \delta\sqrt{n} \Bigr)
\\
&&\qquad= \sum_{j=0}^{n-1}  2  \mathbb{P}\bigl( | N(0,1) | >
\delta\sqrt{n} \bigr)
\nonumber
\\
&&\qquad\le4 n \exp( -\delta^2 n /2 )
\nonumber
,
\end{eqnarray}
where, above, we have used standard Gaussian estimates. Using (\ref
{eqdiff}) and (\ref{equpper}), we finally get
\begin{eqnarray*}
&&\mathbb{P}\bigl( \bigl| \widetilde{H}_{n,k} - \sqrt{p_{\mathrm{max}}}
J_k \bigr| > \varepsilon\bigr)
\\
&&\qquad\le\mathbb{P}\Biggl( 2 \sum_{r=1}^{m-1} \sigma_r \max_{0\le j\le
n-1}  \sup_{0\le t \le1/n} \biggl| \tilde{B}^r\biggl(\frac{j}{n}+t
\biggr) - \tilde{B}^r\biggl(\frac{j}{n}\biggr) \biggr| > \varepsilon\Biggr)
\\
&&\qquad\le\sum_{r=1}^{m-1}  \mathbb{P}\biggl( \max_{0\le j\le n-1}  \sup
_{0\le t \le1/n} \biggl| \tilde{B}^r\biggl(\frac{j}{n}+t\biggr) - \tilde
{B}^r\biggl(\frac{j}{n}\biggr)\biggr | > \frac{\varepsilon}{ 2 \sigma_r
(m-1) } \biggr)
\\
&&\qquad\le4 (m-1) n \exp\biggl( \frac{ -\varepsilon^2 n}{8 \sigma_{\mathrm{max}}^2 (m-1)^2} \biggr)
.
\end{eqnarray*}
\upqed\end{pf}

\section{Uniform binary letters}\label{secbin}

In general, we do not known whether or not the bound in Theorem \ref
{thgen} can be sharpened to $\mathcal{O}(1/\sqrt{n})$. As shown below,
with a more direct proof, for binary alphabets with uniform
distribution this is possible.

Note that for binary alphabets with nonuniform distribution, that is,
for $m=2$ and $k=1$, the limiting distribution $J_1$ is a normal random
variable with zero mean and variance $1- p_{\mathrm{max}}$. Although the
proof of Theorem~\ref{thgen} simplifies in this special case, it still yields
\begin{eqnarray*}
&&\sup_{x\in\mathbb{R}}
\biggl|
\mathbb{P}\biggl( \frac{LI_n - n p_{\mathrm{max}}}{\sqrt{n p_{\mathrm{max}}}} \ge x \biggr)
- \Fi\biggl( \frac{x}{\sqrt{ 1- p_{\mathrm{max}} }} \biggr)
 \biggr|
\\
&&\qquad\le c\biggl ( \sigma_1^2 + \frac{ 1 }{ \sqrt{ p_{\mathrm{max}}
(1-p_{\mathrm{max}}) } } \biggr) \frac{\log n}{\sqrt{n}}
,
\end{eqnarray*}
where $\Fi$ is the standard normal survival function.\vadjust{\goodbreak}

In this section, $m=2$ and assume $\mathbb{P}(X_i=\alpha_1)=\mathbb
{P}(X_i=\alpha_2)=1/2$, $i\in\mathbb{N}$.
Let
\[
Z_i= \cases{
1, & \quad $\mbox{if } X_i=\alpha_1,$
\vspace*{2pt}\cr
-1, & \quad$\mbox{if } X_i=\alpha_2,$}
\]
and let $S_0=0$, $S_k=\sum_{i=1}^k Z_i$, $k\ge1$. Define
\[
\widehat{B}_n(t)= \frac{ S_{[nt]} }{ \sqrt{n} } + (nt-[nt]) \frac{
Z_{[nt]+1} }{ \sqrt{n} }
 ,\qquad  0\le t \le1
.
\]
Then,
\[
\frac{ LI_n - n/2 }{ \sqrt{n} } = -\frac{\widehat{B}_n(1)}{2} + \max
_{t\in[0,1]} \widehat{B}_n(t)
,
\]
and (\ref{eqHtilde}) becomes
\[
\frac{ LI_n - n/2 }{ \sqrt{n} } \Longrightarrow -\frac{B(1)}{2} + \max
_{t\in[0,1]} B(t)
,
\]
where $B$ is a standard Brownian motion.

\begin{theorem}
For any $n\in\mathbb{N}$,
\[
\sup_{x\in\mathbb{R}} \biggl|
\mathbb{P}\biggl( \frac{ LI_n - n/2 }{ \sqrt{n} } \ge x \biggr) -
\mathbb{P}\biggl( -\frac{B(1)}{2} + \max_{t\in[0,1]} B(t) \ge x \biggr)
 \biggr|
\le\frac{24}{\sqrt{n}}
.
\]
\end{theorem}

\begin{pf}
Note that $\max_{t\in[0,1]} \widehat{B}_n(t) = \max_{k=0,\ldots,n}
S_k/\sqrt{n}$. Let
\[
\F(m,b) : = \mathbb{P}\Bigl( \max_{t\in[0,1]} B(t) \ge m, B(1) \le b
\Bigr)
,\qquad  m,b\in\mathbb{R}
\]
and
\begin{eqnarray*}
\F_n(i,j) &:=& \mathbb{P}\Bigl( \max_{k=0,\ldots,n} S_k \ge i, S_n \le
j \Bigr)
\\
&=& \mathbb{P}\biggl( \max_{t\in[0,1]} \widehat{B}_n(t) \ge\frac
{i}{\sqrt{n}}, \widehat{B}_n(1) \le\frac{j}{\sqrt{n}} \biggr)
, \qquad i,j\in\mathbb{Z}
.
\end{eqnarray*}
By the reflection principle, for any $m\ge0$, $b \le m$
\begin{eqnarray*}
\F(m,b) &= &\mathbb{P}\Bigl( \max_{t\in[0,1]} B(t) \ge m, B(1) \ge
m+(m-b) \Bigr)
\\
&=& \mathbb{P}\bigl(B(1) \ge2m-b\bigr) = \Fi(2m-b)
,
\end{eqnarray*}
and for any $i\ge0$, $j \le i$
\begin{eqnarray*}
\F_n(i,j) &= &\mathbb{P}\Bigl( \max_{k=0,\ldots,n} S_k \ge i, S_n \ge
i+(i-j) \Bigr)
\\
&= &\mathbb{P}( S_n \ge2i-j ) = \Fi_n\biggl( 2 \frac{i}{\sqrt
{n}} - \frac{j}{\sqrt{n}} \biggr)
,
\end{eqnarray*}
where
\[
\Fi(z) = \mathbb{P}\bigl(B(1) \ge z\bigr)
,\qquad
\Fi_n(z) = \mathbb{P}\bigl( S_n/\sqrt{n} \ge z \bigr)
,\qquad  z\in\mathbb{R}
.
\]
As is well known (e.g., see~\cite{vanB}),
%
\begin{equation}
\sup_{z\in\mathbb{R}} | \Fi(z)-\Fi_n(z) | \le\frac{ 0.7975
}{ \sqrt{n} }
\label{eqBE}
.
\end{equation}

Next, the joint probability density function of $( \max_{t\in
[0,1]} B(t), B(1) )$ is
\[
f(m,b) = -\frac{ \partial^2 \F(m,b) }{ \partial m\, \partial b } = 2 \Fi''(2m-b)
\]
if $m \ge0$, $b \le m$, and zero elsewhere. For any $x\ge0$, we thus have
\begin{eqnarray}\label{eqcont}
&&\mathbb{P}\biggl( \max_{t\in[0,1]} B(t) -\frac{B(1)}{2} < x \biggr)\nonumber\\
&&\qquad= \int_0^{2x} \int_{2m-2x}^m f(m,b)  \,db\, dm
\nonumber\\
&&\qquad= \int_0^{2x} \int_{2m-2x}^m 2 \Fi''(2m-b)  \,db\, dm
\nonumber
\\[-8pt]
\\[-8pt]
\nonumber
&&\qquad= -2 \int_0^{2x} [  \Fi'(2m-b)  ]_{b=2m-2x}^{b=m}  \,dm\\
&&\qquad= -2 \int_0^{2x} \Fi'(m) - \Fi'(2x)  \,dm
= 2 \Fi(0) -2 \Fi(2x) +2\cdot2x  \Fi'(2x)
\nonumber\\
&&\qquad= 2 \Fi(0) -2 \Fi(2x) -4x  \frac{1}{\sqrt{2\pi}} e^{-{(2x)^2}/{2}}
= 1 -2 \Fi(2x) -4x  \frac{1}{\sqrt{2\pi}} e^{-2x^2}
\nonumber
.
\end{eqnarray}

Observe that $S_n$ is even if $n$ is even, and $S_n$ is odd if $n$ is
odd. In the sequel, \textit{assume that $n$ is even}, in the other case
the computation is similar, and omitted. The joint probability mass
function of $( \max_{k=0,\ldots,n} S_k, S_n )$ is then
\[
p(i,j) = \F_n(i,j) - \F_n(i+1,j) - \F_n(i,j-2) + \F_n(i+1,j-2)
\]
for $j$ even, $i \ge0$, $j \le i$, and zero elsewhere.

For any $x\ge0$, with the notation $l=[x\sqrt{n}]$, we thus have
\begin{eqnarray}\label{eqdisc}
&&\mathbb{P}\biggl( \max_{t\in[0,1]} \widehat{B}_n(t) -\frac{\widehat
{B}_n(1)}{2} < x \biggr)\nonumber\\
&&\qquad= \mathbb{P}\biggl( \max_{k=0,\ldots,n} S_k -\frac{S_n}{2} < l \biggr)
\nonumber\\
&&\qquad= \sum_{i=0}^{2l-2}
\mathop{\sum_{j=2i-2l+2}}_{j \mathrm{even}}^{i}
p(i,j)
\nonumber\\
&&\qquad= \sum_{i=0}^{2l-2} [
\F_n(i,i) - \F_n(i,2i-2l) - \F_n(i+1,i) + \F_n(i+1,2i-2l) ]
\nonumber
\\[-8pt]
\\[-8pt]
\nonumber
&&\qquad= \sum_{i=0}^{2l-2} \biggl[
\Fi_n\biggl( \frac{i}{\sqrt{n}} \biggr) - \Fi_n\biggl( \frac{2l}{\sqrt
{n}} \biggr)
- \Fi_n\biggl( \frac{i+2}{\sqrt{n}} \biggr) + \Fi_n\biggl( \frac
{2l+2}{\sqrt{n}} \biggr) \biggr]\\
&&\qquad= \Fi_n(0) + \Fi_n\biggl( \frac{1}{\sqrt{n}} \biggr)
- \Fi_n\biggl( \frac{2l-1}{\sqrt{n}} \biggr) - \Fi_n\biggl( \frac
{2l}{\sqrt{n}} \biggr)
\nonumber\\
&&\qquad\quad{} -(2l-2) \biggl[ \Fi_n\biggl( \frac{2l}{\sqrt{n}} \biggr) - \Fi
_n\biggl( \frac{2l+2}{\sqrt{n}} \biggr) \biggr]
\nonumber\\
&&\qquad= \Fi_n(0) + \Fi_n\biggl(\frac{2}{\sqrt{n}} \biggr)
- 2 \Fi_n\biggl( \frac{2l}{\sqrt{n}} \biggr)
-(2l-2) \mathbb{P}(S_n = 2l)
\nonumber
;
\end{eqnarray}
where in the last step we used the fact that $\Fi_n$ is constant on the
intervals $[ \frac{i}{\sqrt{n}}, \frac{i+2}{\sqrt{n}} )$,
when $i$ is a nonnegative even integer.

Let us compare (\ref{eqcont}) and (\ref{eqdisc}). Since for any $x\ge
0$, $ 2x\in[ \frac{2l}{\sqrt{n}}, \frac{2l+2}{\sqrt{n}} )$,
by~(\ref{eqBE}),
%
\begin{equation}\label{eqterml}
\sup_{x\ge0}\biggl |  2 \Fi(2x) - 2 \Fi_n\biggl( \frac{2l}{\sqrt{n}}
\biggr)  \biggr|
= \sup_{x\ge0} |  2 \Fi(2x) - 2 \Fi_n(2x)  |
\le\frac{ 1.595 }{ \sqrt{n} }
.
\end{equation}
Moreover, from symmetry considerations, we know that
\[
\Fi_n(0) + \Fi_n\biggl(\frac{2}{\sqrt{n}} \biggr)
= \frac12 + \frac12 \mathbb{P}(S_n = 0) + \frac12 - \mathbb{P}(S_n = 0)
= 1 - \frac12  \mathbb{P}(S_n = 0)
.
\]
Thus,
%
\begin{equation}\label{eqSn0}
1- \biggl( \Fi_n(0) + \Fi_n\biggl(\frac{2}{\sqrt{n}} \biggr) \biggr)
= \frac12  \mathbb{P}(S_n = 0)
= \frac12 \pmatrix{n\vspace*{2pt}\cr n/2}  2^{-n}
.
\end{equation}
Using Stirling's formula
%
\begin{equation}
\sqrt{2\pi}  n^{n+1/2} e^{-n} e^{{1}/{(12n+1)}}
\le n! \le
\sqrt{2\pi}  n^{n+1/2} e^{-n} e^{{1}/{(12n)}}
\label{eqstirling}
,
\end{equation}
the rightmost term in (\ref{eqSn0}) is dominated by
%
\begin{eqnarray}\label{eqterm0}
&&\frac{1}{\sqrt{2\pi}} \sqrt{\frac{n}{ (1/2) n  (1/2) n }}
 \frac{n^n}{( (1/2) n )^{n/2} ( (1/2) n
)^{n/2} }
 e^{{1}/{(12n)} - {1}/{(6n+1}) - {1}/({6n+1})}
 2^{-n}
 \nonumber\hspace*{-35pt}
 \\[-8pt]
 \\[-8pt]
 \nonumber
&&\qquad\le\frac{1}{\sqrt{2\pi}}  \frac{2}{ \sqrt{n} }
\le\frac{0.8}{\sqrt{n}}.\hspace*{-35pt}
\end{eqnarray}
Combining (\ref{eqterml}), (\ref{eqterm0}) and Lemma \ref
{lem3rdterm} below will complete the proof.\vadjust{\goodbreak}
\end{pf}

\begin{lemma}\label{lem3rdterm}
For any $n\in\mathbb{N}$,
%
\begin{equation}\label{eqllt}
\sup_{x\ge0} \biggl|
4x  \frac{1}{\sqrt{2\pi}}  e^{-2x^2} -\bigl (2\bigl[x\sqrt{n}\bigr]-2\bigr) \mathbb
{P}\bigl(S_n = 2\bigl[x\sqrt{n}\bigr]\bigr)
 \biggr|
\le\frac{21}{\sqrt{n}}
.
\end{equation}
\end{lemma}

\begin{pf} First, consider the range $x \ge\sqrt{n}/6$. In this
case, both terms on the left-hand side of (\ref{eqllt}) vanish
exponentially fast as $n\to\infty$. Indeed,
\[
4x  \frac{1}{\sqrt{2\pi}}  e^{-2x^2}
< \frac{4}{\sqrt{2\pi}}  e^{-x^2}
\le\frac{4}{\sqrt{2\pi}}  e^{-n/36}
.
\]
Using $\sqrt{n}  e^{-n/36} \le\sqrt{18}  e^{-18/36} \le\sqrt
{18/e}$, we have
%
\begin{equation}\label{eqexpx}
4x  \frac{1}{\sqrt{2\pi}}  e^{-2x^2}
\le\frac{4  \sqrt{18}}{\sqrt{2\pi e}}  \frac{1}{\sqrt{n}}
= \frac{12}{\sqrt{\pi e}}  \frac{1}{\sqrt{n}}
.
\end{equation}

If $x \ge\sqrt{n}/2 + 1/\sqrt{n}$, $\mathbb{P}(S_n = 2[x\sqrt{n}])=0$.
For $\sqrt{n}/6 \le x < \sqrt{n}/2 + 1/\sqrt{n}$, recalling the
notation $l=[x\sqrt{n}]$, for $n/6 \le l \le n/2$,
\begin{eqnarray}\label{eqchose}
\mathbb{P}(S_n = 2l) &=& \pmatrix{ n\vspace*{2pt}\cr{(n+2l)}/{2}}  2^{-n}
= \pmatrix{n\vspace*{2pt}\cr n/2 + l}  2^{-n}
\nonumber
\\[-8pt]
\\[-8pt]
\nonumber
&\le&\pmatrix{n\vspace*{2pt}\cr n/2 + n/6}  2^{-n}
= \pmatrix{n\vspace*{2pt}\cr 2n/3}  2^{-n}
.
\end{eqnarray}
Using Stirling's formula (\ref{eqstirling}) again, (\ref{eqchose})
can be upper bounded by
\begin{eqnarray*}
&&\frac{1}{\sqrt{2\pi}} \sqrt{\frac{n}{ (2/3) n  (1/3)n }}
 \frac{n^n}{ (( 2/3) n )^{2n/3} ( (1/3) n
)^{n/3} }
 e^{{1}/{(12n)}- {1}/({8n+1}) - {1}/({4n+1})}
 2^{-n}
\\
&&\qquad\le\frac{1}{\sqrt{2\pi}} \frac{9}{ \sqrt{2n} }
 \frac{1}{ [ ( 2/3 )^{2/3} ( 1/3
)^{1/3}  2 ]^n }
 e^{-{123}/{(540n)}}
\le\frac{9}{2 \sqrt{\pi}}  \frac{1}{\sqrt{n}}  e^{-n/18}
.
\end{eqnarray*}
Thus,
\[
(2l-2)  \mathbb{P}(S_n = 2l)
\le\frac{9}{2 \sqrt{\pi}}  \frac{1}{\sqrt{n}}  n e^{-n/18}
.
\]
Since $n e^{-n/18} \le18/e$, we have
%
\begin{equation}
(2l-2)  \mathbb{P}(S_n = 2l)
\le\frac{81}{\sqrt{\pi} e}  \frac{1}{\sqrt{n}}
\label{eqexpl}
.
\end{equation}
Hence, (\ref{eqexpx}) and (\ref{eqexpl}) gives the bound (\ref{eqllt}).

Next, consider the range $0 \le x < \sqrt{n}/6$, with the notation
$l=[x\sqrt{n}]$, $0 \le l < n/6$. The left-hand side of (\ref{eqllt})
can be upper bounded by
\begin{eqnarray}\label{eqsplitted}
&&\biggl|
4x  \frac{1}{\sqrt{2\pi}}  e^{-2x^2} - 4\frac{l}{\sqrt{n}}  \frac
{1}{\sqrt{2\pi}}  e^{-2(l/\sqrt{n})^2}
 \biggr|
\nonumber
\\[-8pt]
\\[-8pt]
\nonumber
&&\qquad{}+ \biggl|
4\frac{l}{\sqrt{n}}  \frac{1}{\sqrt{2\pi}}  e^{-2(l/\sqrt{n})^2} -
(2l-2) \mathbb{P}(S_n = 2l)
 \biggr|
.
\end{eqnarray}
Since the function $xe^{-2x^2}$ is monotone on the intervals $[0,1/2)$
and $[1/2,\infty)$,
\begin{eqnarray}\label{eqabs1}
&&4x  \frac{1}{\sqrt{2\pi}}  e^{-2x^2} - 4\frac{l}{\sqrt{n}}  \frac
{1}{\sqrt{2\pi}}  e^{-2(l/\sqrt{n})^2}
\nonumber\\
&&\qquad\le\frac{4}{\sqrt{2\pi}}  \biggl( \frac{l}{\sqrt{n}} + \frac{1}{\sqrt
{n}} \biggr)   e^{-2( {l}/{\sqrt{n}} + {1}/{\sqrt{n}}
)^2}
- \frac{4}{\sqrt{2\pi}}  \frac{l}{\sqrt{n}}  e^{-2(l/\sqrt{n})^2}
\nonumber\\
&&\qquad\le\frac{4}{\sqrt{2\pi}}  \biggl( \frac{l}{\sqrt{n}} + \frac{1}{\sqrt
{n}} \biggr)   e^{-2( {l}/{\sqrt{n}} + {1}/{\sqrt{n}}
)^2}
- \frac{4}{\sqrt{2\pi}}  \frac{l}{\sqrt{n}}  e^{-2(
{l}/{\sqrt{n}} + {1}/{\sqrt{n}}) ^2}
\\
&&\qquad\le\frac{4}{\sqrt{2\pi}}  \frac{1}{\sqrt{n}}  e^{-2(
{l}/{\sqrt{n}} + {1}/{\sqrt{n}} )^2}
\nonumber\\
&&\qquad\le\frac{4}{\sqrt{2\pi}}  \frac{1}{\sqrt{n}}
\nonumber
.
\end{eqnarray}
On the other hand,
\begin{eqnarray}\label{eqHd}
&&4\frac{l}{\sqrt{n}}  \frac{1}{\sqrt{2\pi}}  e^{-2(l/\sqrt{n})^2} -
4x  \frac{1}{\sqrt{2\pi}}  e^{-2x^2}
\nonumber\\
&&\qquad\le\frac{4}{\sqrt{2\pi}}  \frac{l}{\sqrt{n}}  e^{-2(l/\sqrt{n})^2}
- \frac{4}{\sqrt{2\pi}}  \biggl( \frac{l}{\sqrt{n}} + \frac{1}{\sqrt
{n}} \biggr)   e^{-2( {l}/{\sqrt{n}} + {1}/{\sqrt{n}}
)^2}
\nonumber\\
&&\qquad\le\frac{4}{\sqrt{2\pi}}  \frac{l}{\sqrt{n}}  e^{-2(l/\sqrt{n})^2}
( 1 - e^{- {4l}/{n} - {2}/{n}} )\\
&&\qquad\quad{}- \frac{4}{\sqrt{2\pi}}  \frac{1}{\sqrt{n}}  e^{-2(
{l}/{\sqrt{n}} + {1}/{\sqrt{n}})^2}\nonumber\\
&&\qquad\le\frac{4}{\sqrt{2\pi}}  \frac{l}{\sqrt{n}}  e^{-2l^2/n} ( 1
- e^{- {6l}/{n}} )
- \frac{4}{\sqrt{2\pi}}  \frac{1}{\sqrt{n}}
\nonumber
.
\end{eqnarray}
Using $1-e^{-t} \le t$ $(t\in\mathbb{R})$ with $t=6l/n$, and also $t
e^{-t} \le1/e$ $(t\in\mathbb{R})$ with $t=2l^2/n$, the right-most term
in (\ref{eqHd}) is dominated by
\begin{eqnarray}\label{eqabs2}
&&\frac{12}{\sqrt{2\pi}}  \frac{1}{\sqrt{n}}  \frac{2l^2}{n}  e^{-2l^2/n}
- \frac{4}{\sqrt{2\pi}}  \frac{1}{\sqrt{n}}\nonumber\\
&&\qquad\le\biggl( \frac{12}{e \sqrt{2\pi}}  - \frac{4}{\sqrt{2\pi}}
\biggr) \frac{1}{\sqrt{n}}
\\
&&\qquad\le\biggl( \frac{3}{e}  - 1 \biggr) \frac{4}{\sqrt{2\pi}} \frac
{1}{\sqrt{n}}
\nonumber
.
\end{eqnarray}
From (\ref{eqabs1}) and (\ref{eqabs2}) we get the following bound for
the first term in (\ref{eqsplitted}):
%
\begin{equation}
\biggl|
4x  \frac{1}{\sqrt{2\pi}}  e^{-2x^2} - 4\frac{l}{\sqrt{n}}  \frac
{1}{\sqrt{2\pi}}  e^{-2(l/\sqrt{n})^2}
 \biggr|
\le\frac{4}{\sqrt{2\pi}}  \frac{1}{\sqrt{n}}
\label{eqsplitted1}
.
\end{equation}

To control the second term in (\ref{eqsplitted}), let us recall (see,
e.g., Feller~\cite{Feller}, page~182) that
%
\begin{equation}\label{eqfeller}
\mathbb{P}(S_n = 2l)
= \frac{2}{\sqrt{n}}   \frac{1}{\sqrt{2\pi}}   e^{-{(2l/\sqrt{n})^2}/{2}}
  e^{ \varepsilon_n }
,
\end{equation}
where
%
\begin{equation}\label{eqeps}
-\frac{3l^2}{n^2} - \frac{1}{4n} - \frac{1}{360n^3}
\le\varepsilon_n \le
\frac{2l^4}{n^3} - \frac{1}{4n} + \frac{1}{20n^3}
 \qquad\mbox{if }l<n/6
.
\end{equation}
Hence, for the second term in (\ref{eqsplitted}), we have
\begin{eqnarray}\label{eqsmall}
&&\biggl|
4\frac{l}{\sqrt{n}}  \frac{1}{\sqrt{2\pi}}  e^{-2(l/\sqrt{n})^2} -
(2l-2) \mathbb{P}(S_n = 2l)
 \biggr|
\nonumber\\
&&\qquad\le\frac{2}{\sqrt{2\pi}}  \frac{2l}{\sqrt{n}}  e^{-2(l/\sqrt
{n})^2} |1 - e^{ \varepsilon_n }|
+ 2 \mathbb{P}(S_n = 2l)\\
&&\qquad\le\frac{4}{\sqrt{2\pi}}  \frac{l}{\sqrt{n}}  e^{-2(l/\sqrt{n})^2}
|1 - e^{ \varepsilon_n }|
+ \frac{4}{\sqrt{2\pi}}  \frac{1}{\sqrt{n}}  e^{-2(l/\sqrt{n})^2}
e^{ \varepsilon_n }
\nonumber
.
\end{eqnarray}

If $l\le n^{5/8} /3$, (\ref{eqeps}) becomes
\[
-\frac{1}{3n^{1/3}} - \frac{1}{4n} - \frac{1}{360n^3}
\le\varepsilon_n \le
\frac{2}{81\sqrt{n}} - \frac{1}{4n} + \frac{1}{20n^3}
,
\]
and using $| e^z -1 | \le\max\{ |z|,  | z + \frac
{z^2}{2}  \frac{1}{1-|z|} | \}$, $|z|<1$, (\ref{eqsmall})
can be upper bounded by
\begin{eqnarray}\label{eqeps1}
&&\frac{4}{\sqrt{2\pi}}  \frac{l}{\sqrt{n}}  e^{-2(l/\sqrt{n})^2}
\frac{2}{3 \sqrt{n}}
+ \frac{4}{\sqrt{2\pi}}  \frac{1}{\sqrt{n}}  e^{-2(l/\sqrt{n})^2}
e^{ {2}/{(81 \sqrt{n})} }
\nonumber\\
&&\qquad\le\frac{4}{\sqrt{2\pi}}  \frac{1}{2 \sqrt{e}}  \frac{2}{3 \sqrt{n}}
+ \frac{4}{\sqrt{2\pi}}  \frac{1}{\sqrt{n}}  e^{ {2}/{81} }
\\
&&\qquad\le\frac{4}{\sqrt{2\pi}}  \frac{1.23}{\sqrt{n}}
\nonumber
,
\end{eqnarray}
using $ze^{-2z^2} \le1/(2 \sqrt{e})$ with $z=l/\sqrt{n}$.

If $n^{5/8}/3 < l < n/6$, let us consider (\ref{eqeps}) again and
apply to (\ref{eqsmall}) the trivial upper bound
\begin{eqnarray}\label{eqlast}
&&\frac{4}{\sqrt{2\pi}}  \frac{l}{\sqrt{n}}  e^{-2(l/\sqrt{n})^2} e^{
{2l^4}/{n^3} }
+ \frac{4}{\sqrt{2\pi}}  \frac{1}{\sqrt{n}}  e^{-2(l/\sqrt{n})^2}
e^{ {2l^4}/{n^3} }
\nonumber\\
&&\qquad= \frac{4}{\sqrt{2\pi}}  \frac{l+1}{\sqrt{n}}  e^{ -{2l^2}/{n}
+ {2l^4}/{n^3} }
\\
&&\qquad\le\frac{8}{6 \sqrt{2\pi}}  \frac{1}{\sqrt{n}}  n  e^{ -
{2l^2}/{n} + {2l^4}/{n^3} }
\nonumber
.
\end{eqnarray}
In this range of $l$, it is easy to show that $-\frac{2l^2}{n} + \frac
{2l^4}{n^3} \le-\frac{35}{18} n^{1/4}$, thus (\ref{eqlast}) is itself
dominated by
\[
\frac{8}{6 \sqrt{2\pi}}  \frac{1}{\sqrt{n}}  n  e^{ -{35n^{1/4}}/{18}
 }
.
\]
Using $n e^{ -35n^{1/4}/18 } \le(72/35e)^4$, we further get the upper bounds
%
\begin{equation}\label{eqeps2}
\frac{8}{6 \sqrt{2\pi}} \biggl( \frac{72}{35e} \biggr)^4 \frac
{1}{\sqrt{n}}
\le\frac{0.44}{\sqrt{2\pi}} \frac{1}{\sqrt{n}}
.
\end{equation}

Since (\ref{eqeps1}) is larger than (\ref{eqeps2}), when $0 \le x <
\sqrt{n}/6$, (\ref{eqsplitted1}) and (\ref{eqeps1}) give the
following upper bound for (\ref{eqsplitted}):
\begin{eqnarray*}
&&\biggl|
4x  \frac{1}{\sqrt{2\pi}}  e^{-2x^2} - 4\frac{l}{\sqrt{n}}  \frac
{1}{\sqrt{2\pi}}  e^{-2(l/\sqrt{n})^2}
 \biggr|
\nonumber\\
&&\quad{}+ \biggl|
4\frac{l}{\sqrt{n}}  \frac{1}{\sqrt{2\pi}}  e^{-2(l/\sqrt{n})^2} -
(2l-2) \mathbb{P}(S_n = 2l)
 \biggr|
\nonumber\\
&&\qquad\le\frac{4\cdot2.23}{\sqrt{2\pi}}  \frac{1}{\sqrt{n}}
,
\end{eqnarray*}
which is less than what we had obtained for $x \ge\sqrt{n}/6$.
\end{pf}

%



\printaddresses

\end{document}